\documentclass[11pt]{amsart}

\usepackage{amsmath,amsthm, amscd, amssymb, amsfonts}
\usepackage[all]{xy}
\usepackage[dvips, dvipsnames, usenames]{color}
\allowdisplaybreaks
\usepackage{enumitem}

\newtheorem{theorem}{Theorem}[section]

\newtheorem{lemma}[theorem]{Lemma}
\newtheorem{corollary}[theorem]{Corollary}

\newtheorem{teo intro}{Theorem}

\newtheorem{proposition}[theorem]{Proposition}

\theoremstyle{definition}
\newtheorem{definition}[theorem]{Definition}
\newtheorem{example}[theorem]{Example}

\theoremstyle{remark}
\newtheorem{remark}[theorem]{Remark}

\def\pf{\begin{proof}}
\def\epf{\end{proof}}

\def\zt{\Z^{\theta}}
\def\z2t{\Z^{2\theta}}

\def\yd{{}^{H}_{H}\mathcal{YD}}

\newcommand{\ydg}{{}^{\ku\Gamma}_{\ku\Gamma}\mathcal{YD}}
\newcommand{\ydl}{{}^{\ku\Lambda}_{\ku\Lambda}\mathcal{YD}}

\newcommand{\ydzt}{{}^{\ku\zt}_{\ku\zt}\mathcal{YD}}
\newcommand{\ydztdual}{{}^{(\ku\zt)^*}_{(\ku\zt)^*}\mathcal{YD}}

\newcommand\id{\operatorname{id}}

\newcommand\ord{\operatorname{ord}}

\newcommand\hgt{\operatorname{ht}}

\newcommand\Aut{\operatorname{Aut}}

\newcommand\cop{\operatorname{cop}}

\newcommand\gr{\operatorname{gr}}

\newcommand\op{\operatorname{op}}
\newcommand\ad{\operatorname{ad}}

\newcommand\GKdim{\operatorname{GKdim}}

\def\ot{\otimes}

\def\B{\mathbb{B}}

\newcommand\I{\mathbb I}
\def\L{\mathbb{L}}
\def\N{\mathbb{N}}
\newcommand\Sb{\mathbb S}

\def\Z{\mathbb{Z}}

\def\dpn{\widetilde{\mathcal{B}}}

\def\cB{\mathcal{B}}

\def\cH{\mathcal{H}}
\def\cI{\mathcal{I}}
\def\cJ{\mathcal{J}}

\def\cS{\mathcal{S}}
\def\cU{\mathcal{U}}
\def\cW{\mathcal{W}}

\newcommand{\g}{\mathfrak g}
\def\bq{\mathfrak{q}}

\def\u{\mathfrak{u}}
\def\uti{\widetilde{\mathfrak{u}}}

\def\sfm{\mathsf{m}}
\def\sfM{\mathsf{M}}

\def\sfy{\mathsf{y}}

\newcommand\y[2]{\sfy_{#1}^{(#2)}}

\newcommand{\Eb}{\underline E}
\newcommand{\Fb}{\underline F}
\newcommand{\Kb}{\underline K}
\newcommand{\Lb}{\underline L}

\newcommand{\E}{\mathbf E}

\newcommand{\hb}{\mathbf{h}}
\newcommand{\jb}{\mathbf{j}}
\newcommand{\ku}{ \mathbf{k}}

\def\qmb{\mathbf{q}}
\newcommand{\ub}{\mathbf{u}}
\def\bU{\mathbf{U}}
\newcommand{\vb}{\mathbf{v}}
\def\xb{\mathbf{x}}
\def\yb{\mathbf{y}}
\def\zb{\mathbf{z}}

\newcommand{\base}{\mathtt B}

\newcommand{\Ht}{\mathtt H}

\newcommand{\schi}{{\rho_i(\bq)}}

\newcommand{\pre}{\mathfrak{Pre}}
\newcommand{\post}{\mathfrak{Post}}
\newcommand{\lu}{\mathcal{L}}
\newcommand{\luq}{\lu_{\bq}}

\newcommand{\fO}{\mathfrak O}

\newcommand{\dpndual}{\dpn_{\bq}}
\newcommand{\cBdual}{\cB_{\bq}}

\begin{document}


\title[Quantum divided power algebras]{The quantum divided power algebra of a  finite-dimensional Nichols algebra of diagonal type}

\author[Andruskiewitsch; Angiono; Rossi Bertone]
{Nicol\'as Andruskiewitsch, Iv\'an Angiono, Fiorela Rossi Bertone}

\address{FaMAF-CIEM (CONICET), Universidad Nacional de C\'ordoba,
Medina A\-llen\-de s/n, Ciudad Universitaria, 5000 C\' ordoba, Rep\'
ublica Argentina.} \email{(andrus|angiono|rossib)@mate.uncor.edu}

\thanks{The work was partially supported by CONICET,
FONCyT-ANPCyT, Secyt (UNC)}

\begin{abstract}
Let $\cB_\bq$ be a finite-dimensional Nichols algebra of diagonal type corresponding to a matrix $\bq$.
We consider the graded dual $\lu_{\bq}$ of the distinguished pre-Nichols algebra $\dpn_{\bq}$
from \cite{A-preNichols} and the quantum divided power algebra $\cU_{\bq}$,  a suitable Drinfeld double of $\lu_{\bq} \# \ku \zt$.
We provide basis and presentations by generators and relations of $\lu_{\bq}$ and $\cU_{\bq}$,
and prove that they are noetherian and have finite Gelfand-Kirillov dimension.
\end{abstract}

\maketitle

\section{Introduction}

We fix an  algebraically closed field $\ku$ of characteristic zero.
Let $\g$ be a finite-dimensional simple Lie algebra and $q \in \ku$ a root of 1 (with some restrictions depending on $\g$).
In the theory of quantum groups, there are  several Hopf algebras attached to $\g$ and $q$:
\begin{itemize}[leftmargin=*]\renewcommand{\labelitemi}{$\circ$}
 \item The Frobenius-Lusztig kernel (or small quantum group) $\u_{q}(\g)$.
 \item The $q$-divided power algebra $\cU_{q}(\g)$, see \cite{L-roots of 1, L-libro}.
  \item The quantized enveloping algebra $U_{q}(\g)$, see \cite{DK, DKP, DP}.
\end{itemize}
These Hopf algebras have the following features:
\begin{itemize}[leftmargin=*]\renewcommand{\labelitemi}{$\diamond$}
 \item They admit triangular decompositions, e. g. $\u_{q}(\g) \simeq \u_{q}^+(\g) \otimes \u_{q}^0(\g)\otimes \u_{q}^-(\g)$.
 \item The $0$-part of this triangular decomposition is a Hopf subalgebra, actually a group algebra.
  \item The positive and negative parts are not Hopf subalgebras, but rather Hopf algebras in braided tensor categories, braided Hopf algebras for short.
  \item There are morphisms $\u^+_{q}(\g) \hookrightarrow \cU^+_{q}(\g)$,
  $U_{q}^+(\g) \twoheadrightarrow \u_{q}^+(\g)$  of braided Hopf algebras, and ditto for the full Hopf algebras.
  \item The full Hopf algebras can be reconstructed from the positive part by standard procedures (bosonization, the Drinfeld double).
  \item The  positive part $\u_{q}^+(\g)$ has very special properties-- it is a Nichols algebra.
\end{itemize}

Indeed, $\u_{q}^+(\g)$ is completely determined by the matrix $\bq =(q^{d_ia_{ij}})$, where $(a_{ij})$ is the Cartan matrix of $\g$
and $d_i\in \{1,2,3\}$ make $(d_ia_{ij})$ symmetric. In other words, $\u_{q}^+(\g)$ is the Nichols algebra of diagonal type associated to  $\bq$.

The knowledge of the finite-dimensional Nichols algebras of diagonal type  is crucial in the classification
program of finite-dimensional Hopf algebras \cite{AS Pointed HA}. Two remarkable results on these Nichols algebras are:

\begin{enumerate}\renewcommand{\theenumi}{\alph{enumi}}\renewcommand{\labelenumi}{(\theenumi)}
  \item\label{item:diagonal-classification} The explicit classification  \cite{H-classif RS}.
  \item\label{item:diagonal-relations} The determination of their defining relations
\cite{A-convex, A-presentation}.
\end{enumerate}

Let $\bq\in \ku^{\theta \times \theta}$ with Nichols algebra $\cB_{\bq}$ and assume that $\dim \cB_{\bq} < \infty$.
There are several reasons to consider the analogues of the braided Hopf algebras $U^+_q(\g)$ and $\cU^+_q(\g)$, for $\cB_{\bq}$,
motivated by the classification of Hopf algebras with finite Gelfand-Kirillov dimension and by representation theory.
The analogue $\widetilde\cB_{\bq}$ of  $U^+_q(\g)$
was introduced in \cite{A-presentation} and studied in \cite{A-preNichols} under the name of distinguished pre-Nichols algebra.
The definition of $\widetilde\cB_{\bq}$ is by discarding some of the relations in \cite{A-presentation}.
The purpose of this paper is to study the analogue $\lu_{\bq}$ of $\cU^+_q(\g)$; this is the graded dual
of $\dpndual$ and although it could be called the distinguished post-Nichols algebra of $\bq$,
we prefer to name it the Lusztig algebra as in \cite{AAGTV}, where mentioned in passing.

The paper is organized as follows. Section \ref{sec:preliminaries} is devoted to preliminaries and Section \ref{sec:nichols-diagonal} to
Nichols algebras of diagonal type and distinguished pre-Nichols algebras. In Section \ref{sec:lua} we discuss Lusztig algebras:
we provide a basis and a presentation by generators and relations, and prove that they are noetherian and have finite Gelfand-Kirillov dimension.
In Section \ref{sec:dpa} we introduce the quantum divided power algebra $\cU_{\bq}$, that is a suitable Drinfeld double of $\lu_{\bq} \# \ku \zt$;
we also provide a presentation by generators and relations, and prove that it is noetherian and has finite Gelfand-Kirillov dimension.

\begin{remark} The quantum divided power algebras were introduced and studied in \cite{GH, Hu};
they correspond to Nichols algebras of Cartan type $A_1\times \dots \times A_1$.
\end{remark}

\subsection*{Acknowledgement} We thank the referee for the careful reading of the manuscript.

\section{Preliminaries and conventions}\label{sec:preliminaries}

\subsection{Conventions}\label{subsec:Conventions} If $\theta \in \N$, then we set $\I_{\theta} := \{1, 2, ..., \theta\}$; or simply $\I$ if no confusion arises.
If $\Gamma$ is a group, then $\widehat{\Gamma}$ is its group of characters, that is, one-dimensional representations.

Let $\Sb_{n}$ and $\B_{n}$ be the symmetric and braid groups in $n$ letters, with standard generators $\tau_i = (i\, i+1)$,
respectively $\sigma_i$, $i\in\I_{n-1}$. Let $s: \Sb_{\theta} \rightarrow \B_{\theta}$ be the (Matsumoto) section of the projection
$\pi: \B_{\theta}\twoheadrightarrow \Sb_{\theta}$, $\pi(\sigma_i) = \tau_i$,
$i\in \I_{n-1}$, given by $s(\omega) = \sigma_{i_1}\sigma_{i_2}...\sigma_{i_j}$, whenever
$\omega=\tau_{i_1}\tau_{i_2}...\tau_{i_j}\in\Sb_{\theta}$ has length $j$.

We consider the $\qmb$-numbers in the polynomial ring $\Z[\qmb]$, $ n\in \N$, $0 \leq i \leq n$,
\begin{align*}
(n)_\qmb &=\sum_{j=0}^{n-1}\qmb^{j}, & (n)_\qmb^!&=\prod_{j=1}^{n} (j)_\qmb, &
\binom{n}{i}_\qmb & =\frac{(n)_\qmb^!}{(n-i)_\qmb^!(i)_\qmb^!}.
\end{align*}
If  $q\in\ku$, then $(n)_q$, $(n)_q^!$, $\binom{n}{i}_q$ are the respective evaluations at $q$.

We use the Heynemann-Sweedler notation for coalgebras and comodules; the counit of a coalgebra is denoted by $\varepsilon$,
and the antipode of a Hopf algebra, by $\cS$. All Hopf algebras in this paper have bijective antipode.

Let $H$ be a Hopf algebra.
A \emph{Yetter-Drinfeld module} $V$ over $H$ is a $H$-module and a $H$-comodule satisfying the compatibility condition
\begin{align*}
\delta(h\cdot v) &= h_{(1)}v_{(-1)}\cS (h_{(3)}) \otimes h_{(2)}\cdot v_{(0)}, & h&\in H, v\in V.
\end{align*}
Morphisms of Yetter-Drinfeld modules preserve the action and the coaction.
Thus Yetter Drinfeld modules over $H$ form a braided tensor category $\yd$, with braiding $c_{V,W} (v\otimes w) = v_{(-1)} \cdot w\otimes v_{(0)}$,
$V, W \in \yd$, $v\in V$, $w\in W$. The full subcategory of finite-dimensional objects is rigid.

\subsection{Braided vector spaces and Nichols algebras}\label{subsec:bvs}
A braided vector space is a pair $(V,c)$ where $V$ is a vector space and $c\in \Aut(V\otimes V)$ is a solution of the braid equation
$(c\otimes \id)(\id\otimes c)(c\otimes \id) =
(\id\otimes c)(c\otimes \id)(\id\otimes c)$.

If $V$ is a vector space, then we identify $V^* \otimes V^*$ with a subspace of $(V\otimes V)^*$ by $ \langle f\otimes g,v\otimes w\rangle=
\langle f, w\rangle \langle g, v \rangle$, for $v,w\in V$, $f,g \in V^*$.\footnote{We prefer this identification instead of
$ \langle f\otimes g,v\otimes w\rangle=\langle f, v\rangle \langle g, w \rangle$ because it gives the right extension to tensor categories.}
If $(V,c)$ is a finite-dimensional braided vector space, then $(V^*, c^t)$ is its dual braided vector space, where
$c^t:V^*\otimes V^*\rightarrow V^*\otimes V^*$ is $ \langle c^t(f\otimes g),v\otimes w\rangle=\langle f\otimes g, c(v\otimes w)\rangle$.

We refer to \cite{T} for the basic theory of braided Hopf algebras.
If $R=\bigoplus_{n\geq0} R^n$ is a graded braided Hopf algebra with $\dim R^n < \infty$ for all $n$, then its graded dual
$R^d=\bigoplus_{n\geq0} (R^n)^*$ is again a graded braided Hopf algebra.
We use the variation of the Sweedler notation $\Delta(X)=X^{(1)}\ot X^{(2)}$ for the coproducts in
braided Hopf algebras.

The \emph{Nichols algebra}  of a braided vector space $(V,c)$ is a graded braided Hopf algebra $\cB(V)  = \oplus_{n\ge 0}\cB^n(V)$ with very rigid properties.
There are several alternative definitions of Nichols algebras,
see \cite{AS Pointed HA}. We recall now two of these definitions.

\smallbreak
Let $T(V)= \oplus_{n\ge 0} T^n(V)$ be the  tensor algebra of $V$; it has a braiding $c$ induced from $V$.
Let $T(V) \underline{\otimes} T(V) = T(V) \otimes T(V)$ with the multiplication $(m\otimes m)(\id \otimes c\otimes \id)$ and
let $\Delta: T(V) \to T(V) \underline{\otimes} T(V)$ be the unique algebra map such that $\Delta(v) = v \otimes 1 + 1 \otimes v$, for all $v\in V$.
Then $T(V)$ is a (graded)  braided Hopf algebra with respect to $\Delta$.
Dually, consider the cotensor coalgebra $T^c(V)$ which is isomorphic to $T(V)$ as a vector space.
It bears a multiplication making $T^c(V)$ a braided Hopf algebra with an analogous property, see e. g. \cite{Ro, AG}.
There exists only one morphism of braided Hopf algebras $\Theta:T(V)\rightarrow T^c(V)$ that it is the identity on $V$.
The image of $\Theta$ is the Nichols algebra $\cB(V)$ of $V$.

Here is the second description of $\cB(V)$. Let $\mathfrak{S}$ be the partially ordered set of homogeneous Hopf ideals
of $T(V)$ with trivial intersection with $\ku\oplus V$. Then $\mathfrak{S}$ has a maximal element $\cJ(V)$ and $\cB(V) =  T(V)/\cJ(V)$  \cite{AS Pointed HA}.

\subsection{Pre- and post-Nichols algebras} For several purposes, it is useful to consider braided Hopf algebras $T(V)/I$, for various $I \in \mathfrak{S}$.
These are called \emph{pre-Nichols algebras} \cite{Masuoka}.
Indeed, $\mathfrak{Pre}(V)=\{T(V)/I: I\in\mathfrak{S}\}$ is a poset with ordering given by the surjections; so that it is isomorphic to
$(\mathfrak{S}, \subseteq)$.  The minimal element in $\pre(V)$ is $T(V)$, and the maximal is $\cB(V)$.
Dually, the poset $\post(V)$ consists of graded Hopf subalgebras $S=\bigoplus_{n\geq0} S^n$ of $T^c(V)$ such that $S^1=V$,
ordered by the inclusion. Now the minimal element is $\cB(V)$ and the maximal is $T^c(V)$.
We shall call them \emph{post-Nichols algebras}.

\begin{remark}
The map $\Phi: \pre(V) \rightarrow \post(V^*)$, $\Phi(R) = R^d$, is an anti-isomorphism of posets.
\end{remark}

\begin{proof}
If $R = T(V)/I \in \pre(V)$, then $R^d = I^{\bot}$: hence, $\Phi$ is well-defined and it reverses the order.
Also $\Phi$ is surjective, because for a  given $S\in\post(V^*)$, $I = S^{\bot}$ is a graded Hopf ideal of $T(V)$
and $S = (T(V)/I)^d$.
\end{proof}

\section{Nichols algebras of diagonal type}\label{sec:nichols-diagonal}

A braided vector space $(V,c)$ is of \emph{diagonal type} if there exist a basis $x_1, \dots, x_{\theta}$ of $V$ and a
matrix $\bq = (q_{ij})\in M_{\theta}(\ku^{\times})$ such that $c(x_i\otimes x_j)=q_{ij}x_j\otimes x_i$ for all $i,j \in \I =\I_\theta$.
Let $H = \ku G$ be a group algebra, $\chi_i \in \widehat G$ and $g_j\in Z(G)$ such that $\chi_j(g_i) = q_{ij}$, $i,j \in \I$.
Then $(V, c)$ is realized in $\yd$ by  $h\cdot x_i= \chi_i(h)x_i$ and $\rho(x_i)= g_i\otimes x_i$
for all $i\in \I$, $h\in H$.
We will only consider the case when $H=\ku \zt$,  $g_i = \alpha_i$
and $\chi_j\in\widehat{\zt}$ is given by $\chi_j(\alpha_i)=q_{ij}$, $i, j \in \I$.
Here $\alpha_1, \dots, \alpha_\theta$ is the canonical basis of $\Z^{\theta}$.

Let $V^* \in \ydzt$; it is also a braided vector space of diagonal type, with matrix $\bq$. Indeed, if $y_1, \dots, y_{\theta}$
is the dual basis of $x_1, \dots, x_{\theta}$, then
\begin{align*}
 \langle c^t(y_i\otimes y_j), x_h\otimes x_k\rangle &= \langle y_i\otimes y_j, c(x_h\otimes x_k)\rangle
 = q_{hk} \langle y_i\otimes y_j, x_k\otimes x_h\rangle \\ &=   q_{hk} \delta_{jk}\delta_{ih}
 = q_{ij} \langle y_j\otimes y_i, x_h\otimes x_k\rangle.
\end{align*}

Since $T(V)$ and $\cB_{\bq}=\cB(V)$ are Hopf algebras in $\ydzt$, we may consider the bosonizations $T(V)\# \ku\zt$ and $\cB_{\bq}\# \ku\zt$.
We refer to \cite[\S 1.5]{AS Pointed HA} for the definition of the adjoint action of a Hopf algebra, respectively
the braided adjoint ad$_c$ action of a Hopf algebra in $\ydzt$. Then $\ad_c x\otimes \id = \ad (x\# 1)$ if $x\in T(V)$ or $\cB_{\bq}$,
see \cite[(1-21)]{AS Pointed HA}.

Now the matrix $\bq$ gives rise to a $\Z$-bilinear form $\Xi:\zt\times\zt\to\ku^\times$ by $\Xi(\alpha_j,\alpha_k)=q_{jk}$ for all $j,k\in\I$.
If $\alpha,\beta  \in \zt$, we also set
\begin{align}\label{eq:qab}
 q_{\alpha\beta}  &= \Xi(\alpha,\beta).
\end{align}
The algebra  $T(V)$ is $\zt$-graded. If $x, y \in T(V)$ are homogeneous of degrees $\alpha,\beta  \in \zt$ respectively, then their braided
commutator is
\begin{align}\label{eq:br-comm}
[x,y]_c &= xy - \text{multiplication} \circ c (x\otimes y) = xy -  q_{\alpha\beta} yx.
\end{align}
Note that $\ad_c(x) (y) = [x,y]_c$ whenever $x$ is primitive. We say that  $x$ $\bq$-commutes with a family $(y_i)_{i\in I}$
of homogeneous elements if $[x,y_i]_c=0$, for all $i\in I$.
Same considerations are valid in any braided graded Hopf algebra.

\bigbreak
Define a matrix $(c_{ij}^{\bq})_{i,j\in \I}$ with entries in $\Z \cup \{-\infty \}$ by $c_{ii}^{\bq} = 2$,
\begin{align}\label{eq:defcij}
c_{ij}^{\bq}&:= -\min \left\{ n \in \N_0: (n+1)_{q_{ii}}
(1-q_{ii}^n q_{ij}q_{ji} )=0 \right\},  & i & \neq j.
\end{align}
\emph{We assume from now on that $\dim \cB_{\bq} < \infty$}. Then $c_{ij}^{\bq} \in \Z$ for all $i, j\in \I$
\cite[Section 3.2]{Ro} and we may define the reflections
$s_i^{\bq}\in GL(\Z^\theta)$,  by $s_i^{\bq}(\alpha_j)=\alpha_j-c_{ij}^{\bq}\alpha_i$, $i,j\in \I$.
Let $i\in \I$ and let $\rho_i(V)$ be the braided vector space of diagonal type with matrix  $\rho_i(\bq)$, where
\begin{align}\label{eq:rhoiq}
 \rho_i(\bq)_{jk} &= \Xi(s_i^{\bq}(\alpha_j),s_i^{\bq}(\alpha_k)),& j, k &\in \I.
\end{align}

The proofs of statements (a) and  (b) in the Introduction have as a crucial ingredient the Weyl groupoid
\cite{H-Weyl grp} and the generalized root system \cite{HY};
the definitions involve the assignements $\bq \rightsquigarrow \rho_i(\bq)$ described above.
For our purposes, we just need to recall that
\begin{align}
&\Delta_{\bq}^+\text{ is the set of positive roots of } \cB_{\bq}.
\end{align}

\subsection{Drinfeld doubles}

Let $(V,c)$ be our fixed braided vector space of diagonal type with matrix $\bq$, realized in $\ydzt$ as above.
In this Subsection, the hypothesis on the dimension of the Nichols algebra is not needed.
We describe here the Drinfeld doubles of the bosonizations $T(V)\# \ku\zt$,
$\cB_{\bq}\# \ku\zt$ with respect to suitable bilinear forms. This construction goes back essentially to Drinfeld \cite{Dr-quantum}
and was adapted to different settings in various papers; here we follow \cite{H-isom}.

\begin{definition}
The Drinfeld double $\bU_{\bq}$ of $T(V)\# \ku\zt$ is the algebra generated by elements $E_i$,
$F_i$, $K_i$, $K_i^{-1}$, $L_i$, $L_i^{-1}$, $i\in \I$,
with defining relations
\begin{align*}
XY&=YX, & X,Y \in & \{ K_i^{\pm}, L_i^{\pm}: i\in \I \}, \\
K_iK_i^{-1}&=L_iL_i^{-1}=1,  & E_iF_j-F_jE_i&=\delta_{i,j}(K_i-L_i).
\\ K_iE_j&=q_{ij}E_jK_i,  & L_iE_j&=q_{ji}^{-1}E_jL_i,
\\ K_iF_j&=q_{ij}^{-1}F_jK_i, & L_iF_j&=q_{ji}F_jL_i.
\end{align*}
Then $\bU_{\bq}$ is a $\Z^{\theta}$-graded Hopf algebra, where the comultiplication and the grading are given, for $i\in \I$, by
\begin{align*}
\Delta(K_i^{\pm1})&=K_i^{\pm1} \ot K_i^{\pm1}, & \Delta(E_i)&=E_i \ot 1 + K_i \ot E_i,
\\ \Delta(L_i^{\pm1})&=L_i^{\pm1} \ot L_i^{\pm1}, & \Delta(F_i)&=F_i \ot L_i + 1 \ot F_i.
\\
\deg(K_i)&=\deg(L_i)=0, & \deg(E_i)&=\alpha_i = -\deg(F_i).
\end{align*}
\end{definition}

Let $\bU^+_{\bq}$ (respectively, $\bU^-_{\bq}$) be the subalgebra of $\bU_{\bq}$
generated by $E_i$ (respectively, $F_i$), $i\in \I$. Let $W = (V^*,\bq^t)$. \footnote{Here and in Section \ref{sec:dpa} below, $\bq^t$ corresponds to $V^*$ when realized as Yetter-Drinfeld module over the dual Hopf algebra.}\label{ft:qtrans}
Moreover, $\bU^+_{\bq}$ and $\bU^-_{\bq}$ are Hopf algebras in $\ydzt$ via the actions and coactions
\begin{align*}
  K_i\cdot E_j&=q_{ij}E_j, & &\delta(E_i)=K_i\ot E_i; \\
  L_i\cdot F_j&=q_{ji}F_j, & &\delta(F_i)=L_i\ot F_i.
\end{align*}
Thus, there are isomorphisms $\psi^{+}: T(V) \to \bU^{+}_q$, $\psi^{-}: T(W) \to \bU^{-}_q$ of Hopf algebras in $\ydzt$ given by $\psi^{+}(x_i)=E_i$ and $\psi^{-}(y_i)=F_i$.

Let
\begin{align*}
\u_{\bq}= \bU_{\bq} / (\psi^{-}(\cJ_{\bq^t}) + \psi^{+}(\cJ_{\bq}));
\end{align*}
this is the Drinfeld double of  $\cB_{\bq}\#\ku\Z^\theta$. We denote by $E_i$, $F_i$, $K_i$, $L_i$
the elements of $\u_{\bq}$ that are images of their homonymous in $\bU_{\bq}$.
Let $\u^0$  (respectively, $\u^+_{\bq}$, $\u^-_{\bq}$) be the subalgebra of $\u_{\bq}$ generated by  $K_i$, $L_i$,
(respectively, by  $E_i$, by $F_i$), $i\in \I$. Then $\u^0 \simeq \ku\z2t$;

\begin{itemize}
 \item there is a triangular decomposition $\u_{\bq}\simeq \u^+_{\bq} \otimes \u^0\otimes \u^-_{\bq}$;
 \item  $\u^+_{\bq}\simeq \cB_{\bq}$, $\u^-_{\bq} \simeq \cB_{\bq^t}$.
\end{itemize}

\subsection{Lusztig isomorphisms and PBW bases}

G. Lusztig defined automorphisms of the quantized enveloping algebra $U_q(\g)$ of a simple Lie algebra $\g$, see \cite{L-libro}.
These automorphisms satisfy the relations of the braid group covering the Weyl group of $\g$;
they are instrumental in the construction of Poincar\'e-Birkhoff-Witt (PBW) bases of
$U_q(\g)$. These results were extended to the Drinfeld double of a finite-dimensional Nichols algebra of diagonal type in \cite{H-isom},
with the role of the Weyl group played here by the Weyl groupoid $\cW_{\bq}$. The definition of the
Lusztig isomorphisms in \cite{H-isom} requires some hypotheses on the matrix $\bq$, that are always satisfied in the finite-dimensional case.
So, let $(V,c)$ and $\bq$ as above; recall that we assume that $\dim \cB_{\bq} < \infty$.
Fix $i\in \I$. We first recall the definition of the isomorphisms $\u_{\bq} \to \u_{\schi}$ \cite{H-isom}.
  For $i\neq j\in \I$ and $n\in \N_0$, define the elements of $\u_{\bq}$
\begin{align*}
E_{j,n} &=(\ad E_i)^{n}E_j, &F_{j,n} &=(\ad F_i)^{n}F_j.
\end{align*}

Let $\Eb_j$, $\Fb_j$, $\Kb_j$, $\Lb_j$ be the generators of $\u_{\schi}$. Set
\begin{equation}\label{eqn:escalares lambda}
a_j(\bq):= (-c_{ij}^{\bq})_{q_{ii}}^! \prod_{s=0}^{-c_{ij}^{\bq}-1}(q_{ii}^s q_{ij}q_{ji} - 1) ,  \qquad j \neq i.
\end{equation}

\begin{theorem}{\cite[6.11]{H-isom}}\label{thm: iso Lusztig-Heck}
There are algebra isomorphisms $T_i: \u_{\bq} \to \u_{\schi}$ uniquely determined, for $h, j \in \I$, $j\neq i$, by
\begin{align*}
 T_i(K_h)&=\Kb_i^{-c_{ih}^{\bq}}\Kb_h, &  T_i(E_i)&=\Fb_i\Lb_i^{-1}, & T_i(E_j)&=\Eb_{j,-c_{ij}^{\bq}},
\\
T_i(L_h)&=\Lb_i^{-c_{ih}^{\bq}}\Lb_h, & T_i(F_i)&=\Kb_i^{-1}\Eb_i, & T_i(F_j)&= \frac 1{a_j(\schi)} \, \Fb_{j,-c_{ij}^{\bq}}. \qed
\end{align*}
\end{theorem}

Let $w \in \cW_{\bq}$ be an element of maximal length and  fix a reduced expression $w=\sigma_{i_1}^{\bq} \sigma_{i_2}\cdots \sigma_{i_M}$.
If $k\in \I_M$ and $\hb=(h_1,\dots,h_M)\in\N_0^M$, set
\begin{align} \label{eq:betak}
\beta_k &= s_{i_1}^{\bq}\cdots s_{i_{k-1}}(\alpha_{i_k}),
\\ \label{eq:PBW generators}
E_{\beta_k}&=T_{i_1}\cdots T_{i_{k-1}}(E_{i_k})\in (\u_{\bq}^+)_{\beta_k},
\\ \label{eq:def E a la a, F a la a}
\E^{\hb}&=E_{\beta_M}^{h_M}E_{\beta_{M-1}}^{h_{M-1}} \cdots
E_{\beta_1}^{h_1}.
\end{align}

By \cite[Prop. 2.12]{CH1}, $\Delta_+^{\bq}=\{\beta_k|1\leq k\leq M\}$. Thus, we  set
\begin{align}
\label{eq:Nk} N_\beta &= N_k = \ord q_{\beta\beta} \in\N\cup\{\infty\}, &\text{ if }  \beta &= \beta_k \in \Delta_+^{\bq}.
\end{align}

\begin{theorem} \label{thm: HY PBW bases}
\cite[4.5, 4.8, 4.9]{HY-shapov}
The following set is a basis of $\u^+_{\bq}$:
\begin{align*}
& \{ \E^{\hb}\, | \,  \hb\in\N_0^M, \, 0\leq h_k < N_k, \, k\in \I_M \}. \qed
\end{align*}
\end{theorem}

\subsection{Distinguished pre-Nichols algebra}
We now recall the definition of the distinguished pre-Nichols algebra from \cite{A-preNichols}.
Let $\bq$, $V$ be as above.  First, $i\in\I$ is a \emph{Cartan vertex}  of $\bq$  if
\begin{align}\label{eq:cartan-vertex}
q_{ij}q_{ji} &= q_{ii}^{c_{ij}^{\bq}}, & \text{for all } j \neq i,
\end{align}
recall \eqref{eq:defcij}. Then the set of Cartan roots of $\bq$ is
\begin{align*}
\fO_{\bq} &= \{s_{i_1}^{\bq} s_{i_2} \dots s_{i_k}(\alpha_i) \in \Delta_+^{\bq}:
i\in \I  \text{ is a Cartan vertex of } \rho_{i_k} \dots \rho_{i_2}\rho_{i_1}(\bq) \}.
\end{align*}

A set of defining relations of the Nichols algebra $\cB_{\bq}$, i.~e. generators of the ideal $\cJ_{\bq}$, was given
in \cite[Theorem 3.1]{A-presentation}. We now consider the ideal $\cI_{\bq} \subset \cJ_{\bq}$  of $T(V)$
generated by all the relations in \emph{loc. cit.}, but
\begin{itemize}
 \item we exclude the power root vectors $E_\alpha^{N_\alpha}$, $\alpha\in\fO_{\bq}$,
 \item we add the quantum Serre relations $(\ad_c E_i)^{1-c_{ij}^{\bq}} E_j$ for those $i\neq j$ such that
$q_{ii}^{c_{ij}^{\bq}}=q_{ij}q_{ji}=q_{ii}$.
\end{itemize}

\begin{definition} \cite[3.1]{A-preNichols}
The \emph{distinguished pre-Nichols algebra} of $V$  is
\begin{align*}
\dpn_{\bq}=T(V)/\cI_{\bq}.
\end{align*}
\end{definition}

Let $\uti_{\bq}= \bU_{\bq} / (\psi^{-}(\cI_{\bq^t}) + \psi^{+}(\cI_{\bq}))$; this is
the Drinfeld double of $\dpn_{\bq}\# \ku \zt$. It was shown in \cite{A-preNichols} that
there is a triangular decomposition
$\uti_{\bq}\simeq \uti^+_{\bq} \otimes \uti^0\otimes \uti^-_{\bq}$ as above, with $\uti^0 \simeq \u^0 \simeq \ku\z2t$.

If $\beta_k$ is as in \eqref{eq:betak}, $k\in \I_M$, then we set $\widetilde N_k = \begin{cases} N_k &\mbox{ if }\beta_k\notin\fO_{\bq},
\\ \infty  &\mbox{ if }\beta_k\in\fO_{\bq}, \end{cases}$.
For simplicity, we introduce
\begin{align}\label{eq:ht}
\Ht = \{\hb\in\N_0^M: \, 0\leq h_k < \widetilde N_k, \text{ for all } k\in \I_M \}
\end{align}

\begin{theorem}\label{thm:pnd} \

\begin{enumerate}\renewcommand{\theenumi}{\alph{enumi}}
\renewcommand{\labelenumi}{(\theenumi)}
  \item\label{item:pnd-weyl}  \cite[3.4]{A-preNichols} There exist algebra isomorphisms
 $\widetilde T_i: \uti_{\bq} \to \uti_{\schi}$ inducing  the isomorphisms  $T_i: \u_{\bq} \to \u_{\schi}$.

  \item\label{item:pnd-pbw} \cite[3.6]{A-preNichols} Let $\widetilde{E}_{\beta_k}$, $\widetilde{\E}^{\hb}$ be the elements of $\uti_\bq$ defined as in \eqref{eq:PBW generators}, \eqref{eq:def E a la a, F a la a} with $\widetilde T_i$ instead of $T_i$. Then  
  $\{ \widetilde{\E}^{\hb}\, | \,  \hb\in\Ht\}$ is a basis of $\uti^+_{\bq}$. \qed
\end{enumerate}
\end{theorem}

As before, we have an isomorphism $\widetilde\psi: \dpn_{\bq} \to \uti^{+}_q$ of Hopf algebras in $\ydzt$,
so we define
\begin{align*}
x_{\beta_k} &= \widetilde\psi^{-1}(\widetilde{E}_{\beta_k}),& k&\in \I_M; & \xb^{\hb} &= \widetilde\psi^{-1}(\widetilde{\E}^{\hb}), & \hb &\in \Ht.
\end{align*}
Note that $\widetilde{E}_{\beta_k}$ is a well-defined sequence of braided commutators in the elements $E_i$, $i\in \I$;
then $x_{\beta_k}$ is the same sequence of braided commutators in the $x_i$'s. Also,
$\xb^{\hb} = x_{\beta_M}^{h_M} x_{\beta_{M-1}}^{h_{M-1}} \cdots x_{\beta_1}^{h_1}$ and
$$\base = \{\xb^{\hb}\, | \,  \hb\in\Ht\}$$ is a basis of $\dpn_{\bq}$.
The Hilbert series of a graded vector space $V = \oplus_{n\in \N_0} V^n$ is
$\cH_V =  \sum_{n\in \N_0} \dim V^n \, T^n \in \Z[[T]]$.
It follows  from Theorem \ref{thm:pnd} \eqref{item:pnd-pbw} that
\begin{align}\label{eq:pnd-gk}
\GKdim \dpn_{\bq} &= |\fO_{\bq}|, & \cH_{\dpn_{\bq}} &= \prod_{\beta_k\in\fO_{\bq}} \frac{1}{1 - T^{\deg \beta}}.\prod_{\beta_k\notin\fO_{\bq}} \frac{1 - T^{N_\beta\deg \beta}}{1 - T^{\deg \beta}}.
\end{align}

\section{Lusztig algebras}\label{sec:lua}

Let $\bq = (q_{ij})\in M_{\theta}(\ku^{\times})$, $(V,c)$ the corresponding
braided vector space of diagonal type and $(V^*,\bq)$ the dual braided vector space.
We still assume that $\cB_{\bq}$ is finite-dimensional.
As in \cite[3.3.4]{AAGTV}, we define
the \emph{Lusztig algebra} $\lu_{\bq} $ of $(V,c)$ as the graded dual of the distinguished pre-Nichols algebra $\dpndual$ of $(V^*,\bq)$;
thus, $\cB_{\bq} \subseteq \lu_{\bq}$.
In this Section we establish some basic properties of this algebra.

\subsection{Presentation}\label{subsec:dpa-presentation}
In the rest of the section we consider the bilinear form $\langle \, , \, \rangle: \dpndual \times \dpndual^* \rightarrow \ku$ carried from the identification $V^*\otimes V^*\simeq (V\ot V)^*$ in Section \ref{subsec:bvs} which satisfies for all $x,x' \in \dpndual$, $y,y' \in \dpndual^*$
\begin{align*}
&\langle y, xx'\rangle = \langle y^{(2)},x\rangle\langle y^{(1)},x' \rangle \qquad \mbox{ and }
\qquad \langle yy', x\rangle = \langle y,x^{(2)}\rangle\langle y',x^{(1)} \rangle.
\end{align*}

If $\hb\in\Ht$, then define $\yb_{\hb} \in \dpndual^*$  by $\langle \yb_{\hb},\xb^{\jb}\rangle = \delta_{\hb, \jb}$, $\jb \in \Ht$.
Then $\yb_{\hb} \in \luq$ and $\{ \yb_{\hb}\, | \,  \hb\in\Ht\}$ is a basis of $\luq$.

Let $(\hb_k)_{k\in \I_M}$ denote the canonical basis of $\Z^M$.
If $k\in \I_M$ and $\beta=\beta_k \in \Delta_+^{\bq}$, then we denote the element $\yb_{n\hb_k}$ by $y_{\beta}^{(n)}$.

\smallskip
We recall some notation and results from \cite{A-preNichols} and \cite{AY}.
For $i \in \I_M$, let
\begin{align*}
B^i &= \langle \{ x_{\beta_i}^{h_i} \cdots x_{\beta_1}^{h_1} | 0 \leq h_j < N_j \} \rangle \subseteq \cBdual, \\
\mathbf{B}^i &= \langle \{ x_{\beta_M}^{h_M} \cdots x_{\beta_i}^{h_i} | 0 \leq h_j < N_j \} \rangle \subseteq \cBdual, \\
\widetilde{B}^i &= \langle \{ x_{\beta_i}^{h_i} \cdots x_{\beta_1}^{h_1} | 0 \leq h_j < \widetilde N_j \} \rangle \subseteq \dpndual, \\
\widetilde{\mathbf{B}}^i &= \langle \{ x_{\beta_M}^{h_M} \cdots x_{\beta_i}^{h_i} | 0 \leq h_j < \widetilde N_j \} \rangle \subseteq \dpndual.
\end{align*}

We also denote by $\widetilde{L}^i$ and $\widetilde{\mathbf{L}}^i$ the analogous subspaces of $\luq$:
\begin{align*}
\widetilde{L}^i &= \langle \{ y_{\beta_1}^{(h_1)} \cdots y_{\beta_i}^{(h_i)} | 0 \leq h_j < \widetilde N_j \} \rangle \subseteq \luq, \\
\widetilde{\mathbf{L}}^i &= \langle \{ y_{\beta_i}^{(h_i)} \cdots y_{\beta_M}^{(h_M)} | 0 \leq h_j < \widetilde N_j \} \rangle \subseteq \luq.
\end{align*}

\begin{proposition}\label{prop: coprod cartan}\label{cor: Bi coideal}
\begin{itemize}[leftmargin=*] \item
\cite[4.2, 4.11]{AY}
$B^i$ (respectively $\mathbf{B}^i$) is a right (respectively left) coideal subalgebra of $\cBdual$.

\item \cite[4.1]{A-preNichols}
If $\beta\in\fO_{\bq}$, then $x_{\beta}^{N_{\beta}}$ $\bq$-commutes with every element of $\dpndual$.

\item \cite[4.9]{A-preNichols}
If $\beta_i \in \fO_{\bq}$, then there exist $X(n_1,\dots,n_{i-1})\in \dpndual$ such that
\begin{align*}
\Delta(x_{\beta_i}^{N_{\beta_i}}) &= x_{\beta_i}^{N_{\beta_i}}\ot 1+1\ot x_{\beta_i}^{N_{\beta_i}} \\
&\quad + \sum_{n_k\in \N_0} x_{\beta_{i-1}}^{n_{i-1}N_{\beta_{i-1}}} \dots x_{\beta_1}^{n_1 N_{\beta_1}} \otimes X(n_1,\dots,n_{i-1}).
\hspace{58pt}\qed
\end{align*}

\end{itemize}
\end{proposition}

\begin{corollary} $\widetilde B^i$ is a right coideal subalgebra of $\dpndual$.  \qed
\end{corollary}

Let $Z^+_{\bq}$ be the subalgebra of $\dpndual$ generated by $x_{\beta}^{N_{\beta}}$, $\beta \in \fO_{\bq}$.

\begin{theorem}\cite[4.10, 4.13]{A-preNichols}\label{thm: Z normal subalg}\label{thm: Z Hopf subalg}
$Z^+_{\bq}$ is a braided normal Hopf subalgebra of $\dpndual$. Moreover $Z^+_{\bq} = \,^{co\pi} \dpndual$, where $\pi$ denotes the canonical projection of $\dpndual$ onto $\cB_\bq$.  \qed
\end{theorem}

\begin{lemma}\label{rmk:coprod basis elements}
Let $x$, $x_1$ and $x_2$ be elements in the PBW basis $\base$ of $\dpndual$.
Write $\Delta(x)$ as a linear combination of $\{a\otimes b| a,b \in \base\}$.
Assume that $x_1 \otimes x_2$ has a non-zero coefficient in $\Delta(x)$ (in this combination) and $x_1x_2$ (the concatenation of $x_1$ and $x_2$) is  in $\base$.
Then $x = x_1x_2$.
\end{lemma}
\pf
Suppose that $x=x_{\beta_i}^{h_i}\cdots x_{\beta_1}^{h_1}$ with $h_i>0$. Let
\begin{align*}
 m(x) &= \min\{ j\in\N : h_j\neq 0 \}, \\
 \mathcal{D}(x) &= \sum_{j=1}^{i} \sum_{t=1}^{h_j} \binom{h_j}{t}_{q_{\beta_j\beta_j}} \, x_{\beta_i}^{h_i}\cdots x_{\beta_j}^{t} \otimes x_{\beta_j}^{h_j-t}\cdots x_{\beta_1}^{h_1}
 +1\otimes x, \\
 \widetilde{C}^i &= \langle \{ x_{\beta_M}^{h_M} \cdots x_{\beta_1}^{h_1} \in\base |\, \exists j>i \mbox{ s.t. } h_j\neq0 \} \rangle.
\end{align*}

Observe that if $x_1\ot x_2$ appears in $\mathcal{D}(x)$, then $x=x_1x_2$. However,
if $x_1\ot x_2 \in\sum_{u\in \widetilde B^{i}} u \otimes \widetilde{C}^{m(u)}$, then $x_1x_2\notin\base$.
Therefore the proof is completed by showing that
$$\Delta(x)\in \mathcal{D}(x) + \sum_{u\in \widetilde B^{i}} u \otimes \widetilde{C}^{m(u)}.$$

We proceed by induction on $i$. If $i=1$, then $x=x_{\beta_1}^h$ and $x_{\beta_1}$ is primitive, so
$\Delta(x_{\beta_1}^h)= \sum_{0\le k \le h} \binom{h}{k}_{q_{\beta_1\beta_1}} x_{\beta_1}^k \otimes x_{\beta_1}^{h-k} = \mathcal{D}(x_{\beta_1}^h)$.
Let $i>1$. Now we proceed by induction on $h_i$.
Set $x'=x_{\beta_{i}}^{h_{i}-1}x_{\beta_{i-1}}^{h_{i-1}}\cdots x_{\beta_1}^{h_1}$, so $x=x_{\beta_i}x'$.
Notice that
\begin{align}\label{eq: coproduct e_i}
\Delta(x_{\beta_i})\in  x_{\beta_i}\ot 1 + 1 \ot x_{\beta_i} + \widetilde B^{i-1} \otimes \widetilde{C}^i.
\end{align}
Indeed the analogous statement for $\cBdual$ was proved in \cite[4.3]{AY}, but the same argument applies for $\dpndual$.
By the inductive hypothesis and \eqref{eq: coproduct e_i}
\begin{multline*}
\Delta(x) = \Delta(x_{\beta_i}) \Delta(x') \\
\in \big( x_{\beta_i}\ot 1 + 1 \ot x_{\beta_i} + \widetilde B^{i-1} \otimes \widetilde{C}^i\big)
\Big(\mathcal{D}(x') + \sum_{u\in \widetilde B^{i}} u \otimes \widetilde{C}^{m(u)} \Big) .
\end{multline*}
Notice that $( x_{\beta_i}\ot 1 + 1 \ot x_{\beta_i})\mathcal{D}(x') \in \mathcal{D}(x) + \sum_{u\in \widetilde B^{i}} u \otimes \widetilde{C}^{m(u)}$,
since
\begin{multline*}
(x_{\beta_i}\ot 1 + 1\ot x_{\beta_i})\Big(\sum_{t=1}^{h_i-1} \binom{h_i-1}{t}_{q_{\beta_i\beta_i}} \, x_{\beta_i}^{t} \otimes x_{\beta_i}^{h_i-1-t}\cdots x_{\beta_1}^{h_1} + 1\ot x'\Big)= \\
x_{\beta_i}\ot x' + \sum_{t=2}^{h_i} \binom{h_i-1}{t-1}_{q_{\beta_i\beta_i}} \, x_{\beta_i}^{t} \otimes x_{\beta_i}^{h_i-t}\cdots x_{\beta_1}^{h_1} + \\
 \quad \sum_{t=1}^{h_i-1} q_{\beta_i\beta_i}^{t}\binom{h_i-1}{t}_{q_{\beta_i\beta_i}} \, x_{\beta_i}^{t} \otimes x_{\beta_i}^{h_i-t}\cdots x_{\beta_1}^{h_1} + 1\ot x_{\beta_i}x'
\end{multline*}
and for $h_i>1$, $1\le t< h_i$, we have
$\binom{h_i-1}{t-1}_{q_{\beta_i\beta_i}} + q_{\beta_i\beta_i}^{t}\binom{h_i-1}{t}_{q_{\beta_i\beta_i}}= \binom{h_i}{t}_{q_{\beta_i\beta_i}}$.
Also, $\widetilde B^{i-1}\subset \widetilde B^{i}$, $\widetilde B^{i}$ is a subalgebra and
$\widetilde{C}^iz\subset \widetilde{C}^i$ for all $z\in\dpn_{\bq}$,
by \cite[3.15]{A-preNichols}, so
\begin{align*}
 (\widetilde B^{i-1} \otimes \widetilde{C}^i) \mathcal{D}(x') & \subset \widetilde B^{i-1} \widetilde B^{i} \otimes \widetilde{C}^i \widetilde{B}^i
 \subset \widetilde{B}^i \otimes \widetilde{C}^i.
\end{align*}
As $x_{\beta_i} u \in \widetilde B^{i}$ for all $u\in \widetilde B^{i}$ and $m(u)=m(x_{\beta_i} u)$, then
\begin{align*}
x_{\beta_i} u \otimes \widetilde{C}^{m(u)} =x_{\beta_i} u \otimes \widetilde C^{m(x_{\beta_i} u)} \quad
\mbox{ and } \quad u \otimes  x_{\beta_i}\widetilde{C}^{m(u)} &\subset u \otimes \widetilde C^{m(u)}.
\end{align*}
Finally, $ \widetilde B^{i-1}u  \otimes \widetilde{C}^i \widetilde{C}^{m(u)} \subset \widetilde B^{i} \otimes \widetilde C^i \subset \sum_{v\in \widetilde B^{i}} v \otimes \widetilde{C}^{m(v)}$ for all $u\in \widetilde B^{i}$.
From these considerations the proof of the inductive step follows directly.
\epf

\begin{corollary} If $\beta \in \Delta_+^{\bq}$, then
\begin{align}\label{eq:div power of y}
y_{\beta}^{(r)} &= \frac{y_{\beta}^{r}}{(r)_{q_{\beta\beta}}^!}, & & r<N_{\beta} = \ord q_{\beta\beta}; \\
y_{\beta}^{(n)} &= \frac{(y_{\beta}^{(N_\beta)})^s}{s!} y_{\beta}^{(r)}, & & \beta\in \fO_{\bq}, \, n=sN_\beta + r, \, r<N_\beta.
\end{align}
\end{corollary}
\pf Arguing inductively, we may suppose that
$ y_{\beta}^{r-1} = (r-1)_{q_{\beta\beta}}^! y_{\beta}^{(r-1)}$. If $ x= \xb^{\hb}\in \dpndual$ such that
\begin{align*}
\langle y_{\beta}^{r}, x \rangle =
\langle y_{\beta}^{r-1}, x^{(1)} \rangle \langle y_{\beta}, x^{(2)} \rangle \neq 0,
\end{align*}
then by Lemma \ref{rmk:coprod basis elements}, $x = x_{\beta}^{r}$.
Then
\begin{align*}
\langle y_{\beta}^{r}, x_{\beta}^{r} \rangle =
\langle y_{\beta}^{r-1}, (x_{\beta}^{r})^{(1)} \rangle \langle y_{\beta}, (x_{\beta}^{r})^{(2)} \rangle =
(r-1)_{q_{\beta\beta}}^!(r)_{q_{\beta\beta}} = (r)_{q_{\beta\beta}}^!.
\end{align*}
The second equation follows immediately since $\langle y_{\beta}^{(N_{\beta})}y_{\beta}^{(r)}, x_{\beta}^{N_{\beta}+r} \rangle=1$.
\epf

The next lemma is crucial for the presentation of the algebra $\luq$ by generators and relations.

\begin{lemma}\label{cor: PBW basis luq}
Let $i\in\I_M$, $h_i < \widetilde N_{\beta_i}$ and $\hb=(h_1, \dots, h_M)\in \N_0^M$, then
\begin{align}\label{eq: yb}
\yb_{\hb} &= y_{\beta_1}^{(h_1)}\cdots y_{\beta_M}^{(h_M)}.
\end{align}

Hence
$\{ y_{\beta_1}^{(h_1)} \cdots y_{\beta_M}^{(h_M)} |\, 0 \leq h_i < \widetilde N_{\beta_i} \}$ is a basis of $\luq$.
\end{lemma}
\pf
The proof is by induction on $\operatorname{ht}(\hb) := \sum_{i\in\I_M} h_i$. If $\operatorname{ht}(\hb) =1$ then $\yb_{\hb} = y_\beta$ for some
$\beta \in \Delta_+^{\bq}$ and the claim follows by definition.

Let $1 \leq i_1 < \dots < i_j \leq M$, $n_k < \widetilde N_{\beta_{i_k}}$ and $n_{1}= s N_{\beta_{i_1}} + r \neq 0$ where $r < N_{\beta_{i_1}}$.
Let $y= y_{\beta_{i_1}}^{(n_1)}\dots y_{\beta_{i_j}}^{(n_j)} \in \luq$. Since $\{ \yb_{\hb}\, | \,  \hb\in\Ht\}$ is a basis of $\luq$, we can express $y$ as the linear combination $y = \sum_{\hb\in\Ht} c_{\hb} \yb_{\hb}$. Notice that $c_{\hb}\neq 0$ if and only if $\langle y, x^{\hb} \rangle \neq 0$.

If $r \neq 0$, then we write $y= \frac{1}{(r)_{q}} y_{\beta_{i_1}} y'$ where
$y'= y_{\beta_{i_1}}^{(n_1-1)}\dots y_{\beta_{i_j}}^{(n_j)}$  and $q= q_{\beta_{i_1}\beta_{i_1}}$.
Then $\langle y, x^{\hb} \rangle = \frac{1}{(r)_{q}} \langle y_{\beta_{i_1}}, (x^{\hb})^{(2)} \rangle \langle y', (x^{\hb})^{(1)} \rangle$.
By inductive hypothesis and Lemma \ref{rmk:coprod basis elements}, $c_{\hb}\neq 0$ if and only if $\hb= (0,\dots, n_1, \dots, n_k, 0, \dots)$.
Moreover, the nonzero $c_{\hb}$ is equal to $1$ and the proof in this case is completed.

If $r = 0$, $n_{1} = s N_{\beta_{i_1}}$, then we write $y=y_{\beta_{i_1}}^{(N_{\beta_{i_1}})} y'$.
Arguing as above, \eqref{eq: yb}  follows.
Hence $\{ y_{\beta_1}^{(h_1)} \cdots y_{\beta_M}^{(h_M)} |\, 0 \leq h_i < \widetilde N_{\beta_i} \}$ is a basis of $\luq$
because so is $\{\yb_{\hb} : \hb\in\Ht\}$ by definition.
\epf

We seek for a presentation of $\luq$.
Let us consider the algebra $\L$ presented by generators
$\y{\beta}{n}$, $\beta \in \Delta_+^{\bq}$, $n \in \N$ with relations
\begin{align}\label{eq:dpa-rel-prwnotCartan}
&\y{\beta}{N_\beta} = 0, & \beta &\in \Delta_+^{\bq} - \fO_{\bq};
\\ \label{eq:dpa-rel-prods-power}
&\y{\beta}{h} \y{\beta}{j} = \binom{h+j}{j}_{q_{\beta\beta}} \y{\beta}{h + j}, &
& \begin{matrix} \beta \in \Delta_+^{\bq}, \\ h, j \in \N \end{matrix};
\\ \label{eq:dpa-rel-bracket-power-diff}
&[\y{\beta}{h}, \y{\alpha}{j}]_c = \sum_{\sfm \in \sfM(\alpha, \beta, h, j)} \kappa_{\sfm}\, \sfm, & &
\begin{matrix} \alpha < \beta \in \Delta_+^{\bq},\\0 < h < N_{\alpha},\\ 0 < j  < N_{\beta};\end{matrix}
\\ \label{eq:dpa-rel-bracket-power-diff-cartan}
&[\y{\beta}{N_\beta}, \y{\alpha}{N_\alpha}]_c = \kappa_\gamma \y{\gamma}{N_\gamma} + \hspace{-20pt}
\sum_{\begin{smallmatrix}  0 < l  < N_{\beta}, \, 0 < i  < N_{\alpha} \\  \sfm \in \sfM(\alpha, \beta,N_\alpha-i,N_\beta-l) \end{smallmatrix} } \hspace{-20pt}
\kappa_{\sfm}^{i,l}\, \y{\alpha}{i} \sfm \y{\beta}{l},
& & \begin{matrix} \alpha, \beta, \gamma \in \fO_{\bq}, \\  \alpha <\gamma<\beta;\end{matrix}
\\ \label{eq:dpa-rel-cartan-and-noncartan}
&[\y{\beta}{j}, \y{\alpha}{N_\alpha}]_c =
\sum_{ \begin{smallmatrix}0 < i  < N_{\alpha}, \\ \sfm \in \sfM(\alpha, \beta, N_\alpha-i, j) \end{smallmatrix}}
\kappa_{\sfm}^{i,0}\,\y{\alpha}{i}\sfm , & & \begin{matrix}\alpha \in \fO_{\bq}, \\ \beta \in \Delta_+^{\bq},
\\ 0 < j  < N_{\beta}.\end{matrix}
\end{align}
Here we set
\begin{align*}
\sfM(\alpha, \beta, h, j) &= \{\sfm = \y{\beta_r}{h_r}\cdots \y{\beta_k}{h_k} \in \widetilde L^{\beta} \cap
\widetilde{\mathbf{L}}^{\alpha}: \deg\sfm = \deg \y{\alpha}{h} +\deg \y{\beta}{j} \};
\\
\kappa_{\sfm}^{i,l} &= \langle y_{\beta}^{(h)}y_{\alpha}^{(j)}, x_{\beta}^{l}x_{\beta_k}^{h_k}\cdots x_{\beta_r}^{h_r}x_{\alpha}^{i} \rangle;
\\
\kappa_{\gamma} &= \langle \y{\beta}{N_\beta}\y{\alpha}{N_\alpha}, x_{\gamma}^{N_\gamma} \rangle, \hspace{25pt}
\deg \y{\gamma}{N_\gamma}= \deg \y{\alpha}{N_\alpha}+ \deg \y{\beta}{N_\beta}.
\end{align*}

\begin{theorem}\label{thm:dpa-presentation}
There is an algebra isomorphism $\Upsilon: \L \to \luq$ given by
\begin{align*}
 \Upsilon(\y{\beta}{n}) &= y_{\beta}^{(n)}, & \beta\in \Delta_+^{\bq}, \, n< \widetilde N_\beta.
\end{align*}
\end{theorem}

\pf
We first prove that $\Upsilon$ is well-defined, i.\ e. that \eqref{eq:dpa-rel-prwnotCartan}, \dots,
\eqref{eq:dpa-rel-cartan-and-noncartan} are satisfied by the elements $y_{\beta}^{(n)}\in \luq$.
Relation \eqref{eq:dpa-rel-prwnotCartan} is trivial since $x_{\beta}^{N_\beta} = 0$ if $\beta \notin \fO_{\bq}$
and \eqref{eq:dpa-rel-prods-power} is clear from \eqref{eq:div power of y}.

For the other relations, given $\alpha < \beta$ and $h, j \in\N$, we write
$y_{\beta}^{(h)}y_{\alpha}^{(j)}= \sum_{\hb\in\Ht} c_{\hb} \yb_{\hb}$.
Then $$c_{\hb}= \langle y_{\beta}^{(h)}y_{\alpha}^{(j)}, \xb^{\hb} \rangle
= \langle y_{\alpha}^{(j)}, (\xb^{\hb})^{(1)} \rangle \langle y_{\beta}^{(h)}, (\xb^{\hb})^{(2)} \rangle $$
is the coefficient of $x_{\alpha}^j \ot x_{\beta}^h$ in the expression of $\Delta(\xb^{\hb})$ as linear combination
of elements of the PBW basis in both sides of the tensor product.

If $j<N_\alpha$ and $h<N_\beta$, then $y_{\alpha}^{(j)}, y_{\beta}^{(h)} \in \cB_{\bq}$.
If $c_{\hb} \neq 0$ then $\xb^{\hb}$ appears in the expression of $x_{\alpha}^j x_{\beta}^h$
in elements of the PBW basis, see \cite[Section 3]{A-convex}.
Hence, by \cite[4.8]{HY-shapov} $\xb^{\hb} \in \mathbf{B}^{\alpha} \cap B^{\beta}$,
and relation \eqref{eq:dpa-rel-bracket-power-diff} is clear.

Let $\alpha, \beta\in \fO_{\bq}$, $j=N_{\alpha}$ and $h=N_{\beta}$.
Suppose that there is $\hb = (h_1, \dots, h_M)$ such that $c_{\hb}\neq 0$ and $h_i \geq N_i$ for some $i\in \I_M$.
As $x_{\beta_i}^{N_i}$ $\bq$-commutes with every element of $\dpndual$, we have
$\xb^{\hb}= c\, x_{\beta_i}^{N_i}\xb^{\hb'}$, where $\hb'=(h_1,\dots, h_i-N_i \dots, h_M)$ and
$c= \Xi(h_M\beta_M+\dots +h_{i+1}\beta_{i+1}, N_i\beta_i)\in\ku$. Then
$\Delta(\xb^{\hb}) =  c\,\Delta(x_{\beta_i}^{N_i})\Delta(\xb^{\hb'})$
and hence $\xb^{\hb}=x_{\beta_i}^{N_i}$ by Proposition \ref{prop: coprod cartan}.
For the remaining $\jb$ such that $c_{\jb}\neq 0$ we have $j_i < N_i$ for all $i\in \I_M$. We write
$x_{\alpha}^{N_{\alpha}} \ot x_{\beta}^{N_{\beta}} = \xi (1 \ot x_{\beta}^{n}) (x_{\alpha}^{N_{\alpha}-m} \ot
x_{\beta}^{N_{\beta}-n})(x_{\alpha}^{m} \ot 1)$ where $\xi = \Xi^{-1}((N_{\alpha}-m)\alpha,n\beta)\Xi^{-1}(m\alpha,(N_{\beta}-n)\beta)$.
Therefore, arguing as in the proof of \eqref{eq:dpa-rel-bracket-power-diff} for $y_{\beta}^{(N_{\beta}-n)}y_{\alpha}^{(N_{\alpha}-m)}$,
we obtain that $\yb_{\jb}= \y{\alpha}{m} \sfm \y{\beta}{n}$,
$\sfm \in \widetilde L^{\beta} \cap\widetilde{\mathbf{L}}^{\alpha}$.
Here, either $m=N_\alpha$, $n=N_\beta$ so
$\yb_{\jb}= \Xi(N_{\alpha}\alpha,N_{\beta}\beta)\y{\alpha}{N_\alpha} \y{\beta}{N_\beta}$, or else $m<N_\alpha$ $n<N_\beta$.
Hence relation \eqref{eq:dpa-rel-bracket-power-diff-cartan} follows up to consider the correct degree for $\yb_{\hb}$.

For \eqref{eq:dpa-rel-cartan-and-noncartan}, $ c_{\hb}\neq 0$ implies $\xb^{\hb}\in \cB_{\bq}$ by the same argument above,
since $Z^+_{\bq}$ is a braided Hopf subalgebra by Theorem \ref{thm: Z Hopf subalg}.

Hence, $\Upsilon$ is a morphism of algebras.
By the presentation of $\L$ we can prove that
$\{ \y{\beta_1}{h_1}\dots\y{\beta_M}{h_M} : h_i<\widetilde{N}_i \}$ is a basis of $\L$.
So, $\Upsilon$ maps a basis to a basis by Lemma \ref{cor: PBW basis luq} and then it is bijective.
\epf

\begin{example}
Let $\theta=3 \le N$, $q\in\ku^{\times}$, $\ord q = N$. We consider a diagonal braiding (of super type $A$)
given by a matrix $\bq=(q_{ij})_{i,j \in \I_3}$ such that
\begin{align*}
q_{11}&=q_{23}q_{32}=q, & q_{12}q_{21}&=q^{-1}, & q_{22}&=q_{33}=-1, & q_{13}q_{31}&=1.
\end{align*}
Let $\alpha_{jk} = \displaystyle\sum_{j \le i \le k}\alpha_i$; then
$\Delta^+_\bq = \{\alpha_{jk}:1\le j\le k\le3\}$, $\fO^+_\bq =\{\alpha_1,\alpha_{23},\alpha_{13}\}$.
The Lusztig algebra $\luq$ is presented by generators $y_{jk}^{(n)}$, $1\le j\le k\le3$, $n \in \N$ and relations:

\begin{align*}
y_{12}^{(2)} &=y_2^{(2)}=y_3^{(2)}=0, \\
y_{jk}^{(n)}y_{jk}^{(m)} &= \binom{n+m}{n}_{q_{jk}}y_{jk}^{(n+m)}, \quad n,m \in \N, \\
[y_{12}, y_{1}]_c &= [y_{13}, y_{1}]_c=[y_3, y_1]_c=
[y_{13}, y_{12}]_c =[y_{2}, y_{12}]_c= [y_{23}, y_{12}]_c = 0, \\
[y_{2}, y_{13}]_c &=[y_{23}, y_{13}]_c=[y_{3}, y_{13}]_c=
[y_{23}, y_{2}]_c =[y_{3}, y_{23}]_c= 0, \\
[y_2, y_1]_c &= (1-q^{-1}) y_{12}, \qquad
[y_{3}, y_{12}]_c = (1-q) y_{13}, \\
[y_{23}, y_{1}]_c &= (1-q^{-1}) y_{13}, \qquad
[y_{3}, y_2]_c = (1-q) y_{23}, \\
[y_{23}^{(N)},y_1]_c &= (1-q^{-1})(q_{21}q_{31})^{N-1} y_{13}y_{23}^{(N-1)}, \\
[y_{23},y_{1}^{(N)}]_c &= (1-q^{-1})(q_{21}q_{31})^{N-1} y_1^{(N-1)}y_{13},  \\
[y_2,y_{1}^{(N)}]_c &= (1-q^{-1})q_{21}^{N-1} y_1^{(N-1)}y_{12},\\
[y_{12},y_{1}^{(N)}]_c &= [y_{13},y_1^{(N)}]_c = [y_{3},y_1^{(N)}]_c =0, \\
[y_{13}^{(N)},y_{1}]_c &= [y_{13}^{(N)},y_{12}]_c =[y_{2},y_{13}^{(N)}]_c = [y_{23},y_{13}^{(N)}]_c = [y_{3},y_{13}^{(N)}]_c =0, \\
[y_{23}^{(N)},y_{12}]_c &= [y_{23}^{(N)},y_{13}]_c =[y_{23}^{(N)},y_{2}]_c = [y_{3},y_{23}^{(N)}]_c =0, \\
[y_{13}^{(N)},y_1^{(N)}]_c &= [y_{23}^{(N)},y_{13}^{(N)}]_c = 0, \\
[y_{23}^{(N)},y_1^{(N)}]_c &= (1-q^{-1})^N (q_{21}q_{31})^{N\frac{N-1}{2}} y_{13}^{(N)} \\
&  \qquad + \sum_{k=1}^{N-1} (1-q^{-1})^k (q_{21}q_{31})^{k\frac{2N-k-1}{2}} y_{1}^{(N-k)}y_{13}^{(k)}y_{23}^{(N-k)}.\\
\end{align*}

Indeed, to compute $y_{23}^{(N)}y_1^{(N)}$ in $\luq$, we need to describe all $\hb \in \Ht$, cf. \eqref{eq:ht}, such that
$x_1^{N}\ot x_{23}^{N}$ appears in $\Delta(\xb^\hb)$ with non-zero coefficient (also to be determined), where (for some numeration
of $\Delta^+_\bq$)
\begin{align*}
\xb^{\hb} &= x_3^{h_1}x_{23}^{h_2}x_2^{h_3}x_{123}^{h_4}x_{12}^{h_5}x_1^{h_6}.
\end{align*}
One of these $\xb^\hb$ is $x_{23}^{N}x_1^{N}$, with coefficient $\bq_{N\alpha_1,N\alpha_2+N\alpha_3}$.
Let $\hb$ be as needed. We use the coproduct formulas in \cite[5.1]{A-preNichols}.
Clearly $h_1=0$. From $\Delta(x_{23}^{h_2})$, the only contribution is $(1\ot x_{23})^{h_2}$.
Then we deduce easily that $h_3=h_5=0$, and $h_6=h_2=N-h_4$. In this case, set $h_4=k$
to simplify the notation, so
$$ (1\ot x_{23})^{N-k}(x_1\ot x_{23})^{k}(x_1\ot 1)^{N-k}= (q_{21}q_{31})^{k\frac{2N-k-1}{2}} x_1^{N}\ot x_{23}^{N}. $$
This gives the last relation, and the others are deduced analogously.
\end{example}

\begin{corollary}
The algebra $\luq$ is finitely generated.
\end{corollary}

\pf
By \eqref{eq:dpa-rel-prods-power}, it is generated by  $\{y_{\beta}: \beta\in \Delta_+^{\bq}\} \cup \{y_{\alpha}^{(N_\alpha)}: \alpha\in \fO_{\bq}\}$.
\epf

\begin{remark}
Actually, the subalgebra $\cB_{\bq}\subset\luq$ is generated by its primitive elements
$\{ y_{\alpha}: \alpha \in \Pi_{\bq} \}$ where $\Pi_{\bq}$ denotes the set of simple roots $\alpha_1, \dots, \alpha_{\theta}$.
Moreover, $y_{\gamma}^{(N_\gamma)} \in \ku^{\times} [y_{\beta}^{(N_\beta)}, y_{\alpha}^{(N_\alpha)}]_c$ if and only if $x_{\alpha}^{N_\alpha}\ot x_{\beta}^{N_\beta}$ appears with nonzero coefficient in $\Delta(x_{\gamma}^{N_\gamma})$.
Hence,
\begin{align*}
\{ y_{\alpha}: \alpha\in \Pi_{\bq}\}  \cup  \{y_{\alpha}^{(N_\alpha)}: \alpha\in \fO_{\bq}, \,  x_{\alpha}^{N_{\alpha}}\in\mathcal{P}(\dpndual) \}
\end{align*}
generates $\luq$ as an algebra.
\end{remark}

\begin{proposition}
$\widetilde{\mathbf{L}}^i$ is a right coideal subalgebra of $\luq$.
\end{proposition}
\pf
From Theorem \ref{thm:dpa-presentation} we have that $y_{\beta_j}^{(n)}y_{\beta_i}^{(m)} \in \widetilde{\mathbf{L}}^i$ for $i<j$, thus $\widetilde{\mathbf{L}}^i$ is a subalgebra of $\luq$.
On the other hand, we know that $\langle y_{\beta}^{(n)}, xx'\rangle = \langle (y_{\beta}^{(n)})^{(2)}, x\rangle \langle (y_{\beta}^{(n)})^{(1)}, x'\rangle$. Therefore $y_{\jb}\ot y_{\hb}$ appears with nonzero coefficient in $\Delta(y_{\beta}^{(n)})$ if and only if $x_{\beta}^n$ appears with nonzero coefficient in the expression of $x^{\hb}x^{\jb}$ in the PBW basis. The last condition implies that $x^{\hb}\in \widetilde{B}^{\beta}$ and $x^{\jb}\in\widetilde{\mathbf{B}}^{\beta}$.
Hence,
\begin{align*}
\Delta(y_{\beta}^{(n)}) \in \sum_{i=0}^{n} y_{\beta}^{(i)} \ot y_{\beta}^{(n-i)} +
\widetilde{\mathbf{L}}^{\beta} \ot \widetilde{L}^{\beta}.
\end{align*}
Hence $\Delta(y_{\beta_i}^{(n_i)}\dots y_{\beta_M}^{(n_M)})
= \Delta(y_{\beta_i}^{(n_i)})\Delta(y_{\beta_{i+1}}^{(n_{i+1})}\dots y_{\beta_{M}}^{(n_{M})})\in \widetilde{\mathbf{L}}^i \ot \luq$
and the proof is complete.
\epf

\subsection{Noetherianity and Gelfand-Kirillov dimension}\label{subsec:dpa-noetherian}
We argue as in the pre-Nichols case \cite[Section 3.4]{A-preNichols}, cf. \cite{DP}.
Let us consider the lexicographic order in  $\N_0^{M}$, so that $\hb_M < \dots <\hb_1$, where $(\hb_j)_{j \in \I_M}$
denotes the canonical basis of $\Z^M$.

\begin{lemma}\label{prop: alg filt luq}
Let $\luq(\hb)$ be the subspace of $\luq$ generated by $\yb_{\jb}$, with $\jb \leq \hb$.
Then $\luq(\hb)$ is an $\N_0^M$-algebra filtration of $\luq$.
\end{lemma}
\pf
It is enough to prove that $\yb_{\hb}\yb_{\jb}\in \luq(\hb+\jb)$ for all $\hb, \jb \in \Ht$.
First we consider the case when $\hb = n\hb_k$, $\jb = m\hb_l$, $k, l \in \I_M$, $n,m\in \N$. We claim that $y_{\beta_k}^{(n)}y_{\beta_l}^{(m)} \in \luq(n\hb_k + m\hb_l)$. This follows by definition when $k\leq l$.
If $l<k$, then $[y_{\beta_k}^{(n)},y_{\beta_l}^{(m)}]_c \in \sum_{j<m}y_{\beta_l}^{(j)}.\widetilde{\mathbf{L}}^{l+1}$ by Theorem \ref{thm:dpa-presentation}, thus
\begin{align*}
y_{\beta_k}^{(n)}y_{\beta_l}^{(m)} &\in \luq(n\hb_k + m\hb_l)& &\text{since}& \sum_{j=l+1}^{M} a_j\hb_j &< n\hb_k+m\hb_l.
\end{align*}
The Lemma follows by reordering the factors of $\yb_{\hb}\yb_{\jb}$, for any $\hb, \jb \in \N_0^M$.
\epf

We now consider the corresponding graded algebra
\begin{align*}
\gr\luq &=\oplus_{\hb\in\N_0^M} \gr^{\hb} \luq, & &\text{where}& \gr^{\hb} \luq = \luq(\hb)/\sum_{\jb<\hb} \luq(\jb).
\end{align*}

\begin{lemma}\label{cor: gen rel grluq}
The algebra $\gr\luq$ is presented by generators \emph{$\texttt{y}^{(n)}_k$}, $k \in \I_M$, $n\in\N$, and relations
\emph{\begin{align*}
 &\texttt{y}^{(N_k)}_k=0, \quad \beta_k\notin\fO_{\bq}, \\
 &\texttt{y}^{(n)}_k \texttt{y}^{(m)}_k = \binom{n+m}{m}_{q_{\beta_k\beta_k}} \texttt{y}^{(n+m)}_k, \\
 &[\texttt{y}^{(n)}_k, \texttt{y}^{(m)}_l]_c = 0, \quad l<k.
\end{align*}}
\end{lemma}
\pf
Let $\mathcal{G}$ be the algebra presented by the generators and relations above and
$\pi:\mathcal{G}\rightarrow\gr\luq$ given by $\texttt{y}^{(n)}_k\mapsto y^{(n)}_{\beta_k}$.
By Theorem \ref{thm:dpa-presentation}, the relations above hold in $\gr\luq$. By a direct computation, $\mathcal{G}$
has a basis
\begin{align*}
\{\texttt{y}^{(h_1)}_1\dots\texttt{y}^{(h_M)}_M : h_i<\widetilde{N}_i \}.
\end{align*}
On the other hand, $\yb_{\hb}\in \luq(\hb) - \sum_{\jb<\hb} \luq(\jb)$.
Hence the projection of the PBW basis of $\luq$ is a basis of $\gr\luq$ and $\pi$ is an isomorphism.
\epf

\begin{proposition}
The algebra $\luq$ is Noetherian.
\end{proposition}
\pf
Let $\mathcal{Z}^+$ be the subalgebra of $\gr\luq$ generated by $\{y_{\beta}^{(N_{\beta})}: \beta\in \fO_{\bq}\}$. Then $\mathcal{Z}^+$ is a quantum affine space and $\gr\luq$ is a finitely generated free $\mathcal{Z}^+$-module. Hence $\gr\luq$ is Noetherian and so is $\luq$.
\epf

We compute either from Lemma \ref{cor: PBW basis luq} or else from Lemma \ref{cor: gen rel grluq}  the Gelfand-Kirillov dimension of $\luq$.

\begin{proposition}
 $\GKdim \luq = |\fO_{\bq}|$. \qed
\end{proposition}

\section{Quantum divided power algebras}\label{sec:dpa}

\subsection{Definition}
Let $\bq$, $(V,c)$ be as above with $\dim \cB_{\bq} < \infty$.
Let $W=V^*$, with matrix $\bq^t$, see footnote 2, and let $\{z_{\beta}^{(n)}:\beta\in\Delta_+^{\bq}, n\in\N\}$ be the generators of $\lu_{\bq^t}$.
Here we consider $W\in\ydzt$ via the equivalence of categories between $\ydztdual$ and $\ydzt$.
Then we have a natural evaluation map such that $\langle w \otimes w', v \ot v' \rangle=\langle w\ot v'\rangle \langle w' \ot v \rangle$.
In this section we  define the \emph{quantum divided power algebra} $\cU_{\bq} $ of $(V,c)$
and we  establish some of its basic properties.

Let $\Gamma$ and $\Lambda$ be two copies of $\zt$, generated by $(K_i)_{i\in\I}$ and $(L_i)_{i\in\I}$ respectively; so that
$(K_i^{\pm1})_{i\in\I}$ and $(L_i^{\pm1})_{i\in\I}$ are the generators of $\ku\Gamma$ and $\ku\Lambda$, respectively.
Set $K_{\alpha}= K_1^{a_1}\dots K_{\theta}^{a_{\theta}}$ and $L_{\alpha}=L_1^{a_1}\dots L_{\theta}^{a_{\theta}}$ for
$\alpha=(a_1,\dots,a_{\theta})\in\zt$.
Then $\luq \in\ydg$, $\lu_{\bq^t}\in\ydl$  with structure determined by the formulae
\begin{align*}
&K_{\alpha}^{\pm1}\cdot y_{\beta}^{(n)}= q_{\alpha\beta}^{\pm n} y_{\beta}^{(n)}, &\rho(y_{\beta}^{(n)})= K_{\beta}^n \ot y_{\beta}^{(n)};\\
&L_{\alpha}^{\pm1}\cdot z_{\beta}^{(n)}= q_{\beta\alpha}^{\pm n} z_{\beta}^{(n)},
&\rho(z_{\beta}^{(n)})= L_{\beta}^n \ot y_{\beta}^{(n)}.
\end{align*}
Therefore, we can consider the bosonizations $\luq\#\ku\Gamma$ and $\lu_{\bq^t}\#\ku\Lambda$.

\smallskip
We define next the quantum double of $\luq\#\ku\Gamma$ and $\lu_{\bq^t}\#\ku\Lambda$ following \cite[3.2.2]{Jos}.
For this we need a Hopf pairing between them.

\begin{lemma}
There is a unique bilinear form $(\,|\,):T^c(V) \times (T^c(W))^{\cop} \to \ku$ such that $(1|1)=1$,
\begin{align*}
&(y_i|z_j)=\delta_{ij}, & i,j\in\I; \\
&(y|zz')=(y^{(1)}|z)(y^{(2)}|z'), & y\in T^c(V), \, z,z' \in T^c(W);\\
&(yy'|z)=(y|z^{(1)})(y'|z^{(2)}), & y,y' \in T^c(V), \, z\in T^c(W);\\
&(y|z)=0, &  |y|\neq|z|, \, y\in T^c(V), \, z \in T^c(W).
\end{align*}
\end{lemma}

\pf
Let $\mathbf{T}^n=\sum_{\sigma\in\Sb_n} s(\sigma):(T^c)^n(W)\rightarrow T^n(W)$,
where $s: \Sb_{n} \rightarrow \B_{n}$ is the Matsumoto section, see \cite[\S 3.2]{AG}.
Let $\langle \, , \, \rangle : T^c(V) \ot T(W)^{\op} \rightarrow \ku$ be the evaluation map.
We define $(1|1)=1$,
\begin{align*}
(y|z) &= \langle y , \mathbf{T}^n(z)\rangle, & y \in (T^c)^n(V), z \in (T^c)^n(W)\\
(y|z) &= 0, & y \in (T^c)^n(V), z \in (T^c)^m(W), n\neq m.
\end{align*}

Note that $\mathbf{T}^{i+j}= \mathbf{T}_{i,j}(\mathbf{T}^i\ot\mathbf{T}^j)$ with $\mathbf{T}_{i,j}=\sum s(\sigma^{-1})$
where the sum is over all $(i,j)$-shuffles $\sigma$.
Then, for $y \in (T^c)^n(V)$, $z\in (T^c)^{n-i}(W)$, $z'\in (T^c)^{i}(W)$,
\begin{align*}
(y|zz')& = \langle y , \mathbf{T}^n(z'z)\rangle =
\langle y , \mathbf{T}_{i,n-i}(\mathbf{T}^i\ot\mathbf{T}^{n-i})(z'z)\rangle \\
 & = \langle y , \mathbf{T}_{i,n-i}(\mathbf{T}^i (z') \ot\mathbf{T}^{n-i}(z))\rangle
 =\langle y^{(1)} , \mathbf{T}^{n-i}(z)\rangle \langle y^{(2)} , \mathbf{T}^{i}(z')\rangle \\
 & = (y^{(1)}|z)(y^{(2)}|z')
\end{align*}
The other conditions are clear.
\epf

This bilinear form restricts to $\lu_{\bq} \times (\lu_{\bq^t})^{\cop}$ and then it can be extended to a bilinear form between their bosonizations.
Then we may define a skew-Hopf pairing between $\luq\#\ku\Gamma$ and $\lu_{\bq^t}\#\ku\Lambda$, or equivalently:

\begin{corollary}\label{thm: pairing lusztig alg}
There is a unique Hopf pairing
$$(\,|\,):\luq\#\ku\Gamma \times (\lu_{\bq^t}\#\ku\Lambda)^{\cop} \to \ku$$
such that for all $Y, Y'\in \luq\#\ku\Gamma$,  $Z, Z'\in (\lu_{\bq^t}\#\ku\Lambda)^{\cop}$, $y_{\alpha}^{(n)}\in\luq$, $K_{\alpha}\in\ku\zt$, $z_{\beta}^{(m)}\in\lu_{\bq^t}$ and $L_{\beta}\in\ku\zt$
\begin{align*}
&(Y|ZZ')=(Y_{(1)}|Z)(Y_{(2)}|Z'), \quad
(YY'|z)=(Y|Z_{(1)})(Y'|Z_{(2)}), \\
&(y_{\alpha}^{(n)}|z_{\beta}^{(m)})= \delta_{n\alpha,m\beta}, \quad (y_{\alpha}^{(n)}|L_{\beta})=0,  \quad
(K_{\alpha}|z_{\beta}^{(m)})=0, \quad
(K_{\alpha}|L_{\beta})= q_{\alpha\beta}.
\end{align*}
Moreover, this pairing satisfies the equation $(yK|zL)=(y|z)(K|L)$. \qed
\end{corollary}

Let $\cU_{\bq}$ be
the Drinfeld double of $\lu_{\bq}\#\ku\Gamma$ and $(\lu_{\bq^t}\#\ku\Lambda)^{\cop}$
with respect to the Hopf pairing in Corollary \ref{thm: pairing lusztig alg}.
In other words:

\begin{definition}
Let $\cU_{\bq}$ be the unique Hopf algebra such that
\begin{enumerate}
\item $\cU_{\bq}= (\lu_{\bq}\#\ku \Gamma) \ot (\lu_{\bq^t}\#\ku \Lambda)$ as vector spaces,
\item the maps $Y \mapsto Y\ot 1$ and $Z\mapsto 1\ot Z$ are Hopf algebra morphisms,
\item the product is given by
$$(Y\ot Z)(Y'\ot Z')= (Y'_{(1)}|\cS(Z_{(1)}))YY'_{(2)} \ot Z_{(2)}Z'(Y'_{(3)}|Z_{(3)}) $$ for all $Y,Y'\in\luq\#\ku\Gamma$ and $Z,Z'\in(\lu_{\bq^t}\#\ku\Lambda)^{\cop}$.
\end{enumerate}
\end{definition}

By the construction of $\cU_{\bq}$, there is a triangular decomposition, via the multiplication,
$\cU_{\bq}\simeq \cU^+_{\bq} \otimes \cU^0\otimes \cU^-_{\bq}$ where
\begin{align*}
\cU^+_{\bq}&\simeq \lu_{\bq},& \cU^-_{\bq} &\simeq \lu_{\bq^t},& \cU^0 \simeq \ku(\zt \times \zt).
\end{align*}

\smallskip
We give a presentation of the algebra $\cU_{\bq}$ by generators and relations.
The tensor product signs in elements of $\cU_{\bq}$ will be omitted.

\begin{proposition}\label{prop:gen rel Uq}
The algebra $\cU_{\bq}$ is generated by the elements $y_{\beta}^{(n)}$, $z_{\beta}^{(n)}$, $K^{\pm 1}_\beta$,
$L^{\pm 1}_\beta$ for $\beta \in \Delta^{\bq}_+$, $n\in\N$; and relations
\eqref{eq:dpa-rel-prwnotCartan}, \dots, \eqref{eq:dpa-rel-cartan-and-noncartan} between the $y_{\beta}^{(n)}$'s,
similar relations for the $z_{\beta}^{(n)}$'s plus the relations
\begin{align}
K_\beta K^{-1}_\beta = L^{-1}_\beta L_\beta= 1, \qquad
K^{\pm 1}_\beta L^{\pm 1}_\alpha = L^{\pm 1}_\alpha K^{\pm 1}_\beta \\
K_\alpha y_{\beta}^{(n)} = q_{\alpha\beta}^n y_{\beta}^{(n)} K_\alpha, \qquad
L_\alpha y_{\beta}^{(n)}= q_{\beta\alpha}^{-n} y_{\beta}^{(n)} L_\alpha, \\
K_\alpha z_{\beta}^{(n)} = q_{\alpha\beta}^{-n} z_{\beta}^{(n)} K_\alpha, \qquad
L_\alpha z_{\beta}^{(n)}= q_{\beta\alpha}^n z_{\beta}^{(n)} L_\alpha, \\
zy = ({y}^{(1)}|\cS ({z}^{(3)})) \, (K_2K_3|L_3^{-1}) \, ({y}^{(3)}|{z}^{(1)}) \,
{y}^{(2)} K_{3} {z}^{(2)} L_{3} \label{eq:rel Uq},
\end{align}
for all $\alpha,\beta \in \Delta^{\bq}_+$, $n,m\in\N$.
Here in \eqref{eq:rel Uq} $y= y_{\beta}^{(n)}\in\lu_{\bq}$, $z=z_{\alpha}^{(m)}\in\lu_{\bq^t}$, and denote
$K_{i}= ({y}^{(i)})_{(-1)}$ and $L_{i}=({z}^{(i)})_{(-1)}$ for the coactions of $\ku\Gamma$ and $\ku\Lambda$ respectively.
\qed\end{proposition}

Note that if $y= y_{\alpha_i}$, $z= z_{\alpha_j}$ with $\alpha_i, \alpha_j \in \Pi_{\bq}$,
then $y$, $z$ are primitives and relation \eqref{eq:rel Uq} is
$zy - yz = \delta_{ij} (K_i - L_i)$.

\subsection{Basic properties}

Proceeding as in \cite{DP, A-preNichols}, we will prove that the algebra $\cU_{\bq}$ is Noetherian.
For each $\hb, \jb \in \Ht$, $K\in \Gamma$, $L\in \Lambda$, set
\begin{align*}
d_1(\yb_{\hb}KL\zb_{\jb}) &= \sum_{i\in\I_M}(h_i+j_i)\hgt(\beta_i), \\
d(\yb_{\hb}KL\zb_{\jb}) &= \Big(d_1(\yb_{\hb}KL\zb_{\jb}),\, h_1, \dots, h_M, \, j_1, \dots, j_M \Big) \in\N_0^{2M+1}.
\end{align*}

Consider the lexicographic order in $\N_0^{2M+1}$. If $\ub \in \N_0^{2M+1}$, then we set
\begin{align*}
\cU_{\bq}(\ub) &= \text{span of } \{\yb_{\hb}KL\zb_{\jb}: \hb, \jb \in \Ht,\, K\in \Gamma, \, L\in \Lambda, \,
d(\yb_{\hb}KL\zb_{\jb}) \leq \ub \}.
\end{align*}

\begin{lemma}
$(\cU_{\bq}(\ub))_{\ub \in \N_0^{2M+1}}$ is an $\N_0^{2M+1}$-algebra filtration of $\cU_{\bq}$.
\end{lemma}
\noindent\emph{Proof.}
It is enough to prove that $(\yb_{\hb}KL\zb_{\jb})(\yb_{\hb'}K'L'\zb_{\jb'})\in \cU_{\bq}(\ub+\ub')$ for all $\hb, \jb, \hb', \jb' \in \Ht$, $K, K'\in\Gamma$ and $L, L'\in\Lambda$ where $d(\yb_{\hb}KL\zb_{\jb})=\ub$ and $d(\yb_{\hb'}K'L'\zb_{\jb'})=\ub'$.

First we claim that
\begin{equation}\label{eq:ht(zy-yz)}
d_1(z_{\beta}^{(n)}y_{\alpha}^{(m)} - y_{\alpha}^{(m)}z_{\beta}^{(n)}) < m\hgt(\alpha) + n\hgt(\beta).
\end{equation}
Indeed, since the coproduct in $\luq$ (resp. $\lu_{\bq^t}$) is graded, we have that
$d_1((y_{\alpha}^{(m)})^{(2)}) < m\hgt(\alpha)$ if $(y_{\alpha}^{(m)})^{(1)}\neq 1$
(resp. $d_1((z_{\beta}^{(n)})^{(2)}) < n\hgt(\beta)$ if $(z_{\beta}^{(n)})^{(1)}\neq 1$).
Hence, for $K\in\Gamma$ and $L\in\Lambda$ we have
$$d_1((y_{\alpha}^{(m)})^{(2)}KL(z_{\beta}^{(n)})^{(2)})\leq m\hgt(\alpha) + n\hgt(\beta)$$ and by Proposition
\ref{prop:gen rel Uq} the claim follows.

Since $K$, $L$ $\bq$-commutes with all elements of $\luq$ and $\lu_{\bq^t}$ for all $K\in\Gamma$ and $L\in\Lambda$.
We proceed as in Lemma \ref{prop: alg filt luq} and we reduce the proof to the product between $z_{\beta_i}^{(n)}$ and $y_{\beta_j}^{(m)}$.
It follows directly by \eqref{eq:ht(zy-yz)} that
$$z_{\beta_i}^{(n)}y_{\beta_j}^{(m)}\in \cU_{\bq}(m\hgt(\beta_j) + n\hgt(\beta_i), \delta_j,\delta_i). \qed$$

We consider the associated graded algebra $\gr\cU_{\bq} =\oplus_{\vb\in\N_0^{2M+1}} {\cU_{\bq}}^{\vb}$
where ${\cU_{\bq}}^{\vb}= \cU_{\bq}(\vb)/\sum_{\ub<\vb} \cU_{\bq}(\ub)$.

\begin{corollary}\label{cor: gen rel grUq}
The algebra $\gr\cU_{\bq}$ is presented by generators \emph{$\texttt{y}^{(n)}_j$, $\texttt{z}^{(n)}_j$, $K^{\pm 1}_j$, $L^{\pm 1}_j$},
$j \in \I_M$, $n\in\N$ and relations

\emph{\begin{align*}
 & RS=SR, & R,S\in\{K^{\pm 1}_j,L^{\pm 1}_j : j\in\I_M  \}\\
 & K_\beta K^{-1}_\beta = L_\beta L^{-1}_\beta= 1
 & \texttt{y}^{(n)}_k \texttt{z}^{(m)}_l = \texttt{z}^{(m)}_l \texttt{y}^{(n)}_k\\
 &\texttt{y}^{(N_k)}_k=0, \quad \beta_k\notin\fO_{\bq}, &\texttt{z}^{(N_k)}_k=0, \quad \beta_k\notin\fO_{\bq}, \\
 &\texttt{y}^{(n)}_k \texttt{y}^{(m)}_k = \binom{n+m}{m}_{q_{\beta_k\beta_k}} \texttt{y}^{(n+m)}_k,
 &\texttt{z}^{(n)}_k \texttt{z}^{(m)}_k = \binom{n+m}{m}_{q_{\beta_k\beta_k}} \texttt{z}^{(n+m)}_k,\\
 &[\texttt{y}^{(n)}_k, \texttt{y}^{(m)}_l]_c = 0, \quad l<k,
 &[\texttt{z}^{(n)}_k, \texttt{z}^{(m)}_l]_c = 0, \quad l<k, \\
 & K_\alpha y_{\beta}^{(n)}= q_{\alpha\beta}^n y_{\beta}^{(n)} K_\alpha,
 & K_\alpha z_{\beta}^{(n)}= q_{\alpha\beta}^{-n} z_{\beta}^{(n)} K_\alpha, \\
 & L_\alpha y_{\beta}^{(n)}= q_{\beta\alpha}^{-n} y_{\beta}^{(n)} L_\alpha,
 & L_\alpha z_{\beta}^{(n)}= q_{\beta\alpha}^n z_{\beta}^{(n)} L_\alpha.
  \end{align*}}
\end{corollary}
\pf
The proof of this statement is similar to the proof of Lemma \ref{cor: gen rel grluq} if we check
that $y^{(n)}_k z^{(m)}_l = {z}^{(m)}_l {y}^{(n)}_k$ for all $y^{(n)}_k\in\luq$ and $z^{(m)}_l\in\lu_{\bq^t}$;
but this follows by \eqref{eq:ht(zy-yz)}.
\epf

\begin{proposition}
The algebra $\cU_{\bq}$ is Noetherian and  $\GKdim \cU_{\bq} = 2|\fO_{\bq}|+2\theta$.
\end{proposition}
\noindent\emph{Proof.}
Let $\mathcal{Z}$ be the subalgebra of $\gr\cU_{\bq}$ generated by
$\{K_i, L_i: i\in\I\}$ and $\{ y_{\beta}^{(N_{\beta})}, z_{\beta}^{(N_{\beta})}: \beta\in \fO_{\bq}\}$. Then $\mathcal{Z}$ is the localization of a quantum affine space and $\gr\cU_{\bq}$ is a free $\mathcal{Z}$-module of rank $\prod_{i\in\I_M} N_i$.
Therefore $\gr\cU_{\bq}$ is Noetherian and so is $\cU_{\bq}$.
Moreover, by \cite[Prop. 6.6]{KL},
$$\GKdim \cU_{\bq} = \GKdim \gr\cU_{\bq}= \GKdim \mathcal{Z} = 2|\fO_{\bq}|+2\theta.\qed$$

\end{document}